\documentclass[reqno,10pt,centertags,draft]{amsart}
\usepackage{amsmath,amsthm,amscd,amssymb,latexsym,upref,enumerate}
\newcommand{\bbC}{{\mathbb{C}}}

\newcommand{\bbN}{{\mathbb{N}}}
\newcommand{\bbR}{{\mathbb{R}}}

\newcommand{\bbZ}{{\mathbb{Z}}}

\newcommand{\cC}{{\mathcal C}}

\newcommand{\no}{\notag}
\newcommand{\lb}{\label}
\newcommand{\f}{\frac}

\newcommand{\ol}{\overline}

\newcommand{\wti}{\widetilde}

\newcommand{\oh}{o}

\newcommand{\supp}{\text{\rm{supp}}}

\newcommand{\bi}{\bibitem}

\newcommand{\deven}{\delta_{\rm even}}
\newcommand{\dodd}{\delta_{\rm odd}}

\renewcommand{\Re}{\text{\rm Re}}
\renewcommand{\Im}{\text{\rm Im}}

\newcommand{\abs}[1]{\lvert#1\rvert}

\newcommand{\si}{\sigma}

\newcommand{\al}{\alpha}

\newcommand{\de}{\delta}

\newcommand{\Te}{\Theta}
\newcommand{\ze}{\zeta}

\newcommand{\veps}{\varepsilon}

\newcommand{\C}{\mathbb{C}}
\newcommand{\Cz}{{\C\backslash\{0\}}}
\newcommand{\CdD}{{\C\backslash\dD}}

\newcommand{\D}{\mathbb{D}}
\newcommand{\dD}{{\partial\hspace*{.2mm}\mathbb{D}}}
\newcommand{\Z}{{\mathbb{Z}}}
\newcommand{\N}{{\mathbb{N}}}

\newcommand{\Lt}[1]{{L^2(\dD;d\mu_{#1}(\cdot,k_0))}}

\newtheorem{theorem}{Theorem}[section]

\newtheorem{lemma}[theorem]{Lemma}
\newtheorem{corollary}[theorem]{Corollary}
\newtheorem{hypothesis}[theorem]{Hypothesis}
\theoremstyle{definition}

\allowdisplaybreaks \numberwithin{equation}{section}

\begin{document}

\title[Borg--Marchenko-type Uniqueness results for CMV operators]
{Borg--Marchenko-type Uniqueness Results for CMV Operators}

\author[S.\ Clark, F.\ Gesztesy, and M.\ Zinchenko]
{Stephen Clark, Fritz Gesztesy, and Maxim Zinchenko}

\address{Department of Mathematics \& Statistics,
University of Missouri, Rolla, MO 65409, USA}
\email{sclark@umr.edu}
\urladdr{http://web.umr.edu/\~{}sclark/index.html}

\address{Department of Mathematics,
University of Missouri, Columbia, MO 65211, USA}
\email{fritz@math.missouri.edu}
\urladdr{http://www.math.missouri.edu/personnel/faculty/gesztesyf.html}

\address{Department of Mathematics,
California Institute of Technology, Pasadena, CA 91125, USA}
\email{maxim@caltech.edu}
\urladdr{http://www.math.caltech.edu/\~{}maxim}

%%%%%%%%%%%%%%%%%%%%%%%%%%%%%%%%%%%%%%%%
\dedicatory{Dedicated with great pleasure to Pavel Exner on the occasion 
of his 60th birthday}
%%%%%%%%%%%%%%%%%%%%%%%%%%%%%%%%%%%%%%%%

\thanks{Based upon work supported by the US National Science
Foundation under Grants No.\ DMS-0405526 and DMS-0405528.}
\subjclass{Primary 34E05, 34B20, 34L40;  Secondary 34A55}
\keywords{CMV operators, orthogonal polynomials, finite difference operators, Weyl--Titchmarsh theory, Borg--Marchenko-type uniqueness theorems.}
%\thanks{\it .}

%\date{\today}
%\date{March 19, 2007.}

%%%%%%%%%%%%%%%%%%%%%%%%%%%%%%%%%%%%%%%%
\begin{abstract}
We prove local and global versions of Borg--Marchenko-type uniqueness theorems for 
half-lattice and full-lattice CMV operators (CMV for Cantero, Moral, and Vel\'azquez 
\cite{CMV03}). While our half-lattice results are formulated 
in terms of Weyl--Titchmarsh functions, our full-lattice results involve the diagonal and main off-diagonal Green's functions.
\end{abstract}
%%%%%%%%%%%%%%%%%%%%%%%%%%%%%%%%%%%%%%%%

\maketitle

%%%%%%%%%%%%%%%%%%%%%%%%%%%%%%%%%%%%%%%%
%%%%%%%%%%%%%%%%%%%%%%%%%%%%%%%%%%%%%%%%
\section{Introduction}\lb{s1}
%%%%%%%%%%%%%%%%%%%%%%%%%%%%%%%%%%%%%%%%
%%%%%%%%%%%%%%%%%%%%%%%%%%%%%%%%%%%%%%%%

To set the stage, we briefly review the history of Borg--Marchenko-type uniqueness theorems. Apparently, it all started in connection with Schr\"odinger operators on half-lines, and so we turn to that case first.

Let $H_j = -\f{d^2}{dx^2} + V_j$, $V_j\in L^1 ([0,R]; dx)$ for all $R>0$, $V_j$ 
real-valued, $j=1,2$, be two self-adjoint operators in $L^2 ([0,\infty); dx)$ which, just for simplicity, have a Dirichlet boundary condition at $x=0$ (and possibly a self-adjoint boundary condition at infinity). Let $m_j(z)$, $z\in\bbC\backslash\bbR$, be the Weyl-Titchmarsh $m$-functions associated with $H_j$, $j=1,2$. Then the celebrated Borg--Marchenko uniqueness theorem, in this particular context, reads as follows:

%%%%%%%%%%%%%%%%%%%%%%%%%%%%%%%%%%%%%
\begin{theorem} \lb{t1.1} %\mbox{\rm (\cite{Bo52,Ma50,Ma52})} 
Suppose 
\begin{equation}  
m_1(z) = m_2(z), \;\; z\in\bbC\backslash\bbR, \, \text{ then } \, 
V_1(x) = V_2(x) \, \text{ for a.e.\  $x\in [0,\infty)$.}   \lb{1.1}
\end{equation}
\end{theorem}
%%%%%%%%%%%%%%%%%%%%%%%%%%%%%%%%%%%%%

This result was published by Marchenko \cite{Ma50} in 1950. Marchenko's extensive treatise on spectral theory of one-dimensional Schr\"odinger operators \cite{Ma52}, repeating the proof of his uniqueness theorem, then appeared in 1952, which also marked the appearance of Borg's proof of the uniqueness theorem \cite{Bo52} (apparently, based on his lecture at the 11th Scandinavian Congress of Mathematicians held at Trondheim, Norway in 1949).

As pointed out by Levitan \cite{Le87} in the Notes to Chapter~2, Borg and Marchenko  were actually preceded by Tikhonov \cite{Ti49} in 1949, who proved a special case of  Theorem~\ref{t1.1} in connection with the string equation (and hence under certain  additional hypotheses on $V_j$). Since Weyl--Titchmarsh functions $m(z)$ are uniquely  
related to the spectral measure $d\rho$ of a self-adjoint (Dirichlet) Schr\"odinger 
operator $H=-\f{d^2}{dx^2} + V$ in $L^2 ([0,\infty))$ by the standard Herglotz  representation
\begin{equation} \lb{1.2}
m(z) = \Re (m(i)) + \int_\bbR d\rho(\lambda) 
[(\lambda -z)^{-1} - 
\lambda(1+\lambda^2)^{-1}], \quad z\in\bbC\backslash\bbR,
\end{equation}
Theorem \ref{t1.2} is equivalent to the following statement: 
Denote by $d\rho_j$ the 
spectral measures of $H_j$, $j=1,2$. Then
\begin{equation} \lb{1.3}
d\rho_1 = d\rho_2 \, \text{ implies $V_1 = V_2$ a.e.\ on $[0,\infty)$.}
\end{equation}
In fact, Marchenko's proof takes the spectral measures $d\rho_j$ as the point of 
departure while Borg focuses on the Weyl--Titchmarsh functions $m_j$.

We emphasize at this point that Borg and Marchenko also treat the general case of 
non-Dirichlet boundary conditions at $x=0$ (in which equality of the two $m$-functions also identifies the two boundary conditions); moreover, Marchenko also simultaneously discussed the half-line and the finite interval case. For brevity we chose to illustrate the simplest possible case only.

To the best of our knowledge, the only alternative approaches to Theorem \ref{t1.1} 
are based on the Gelfand--Levitan solution \cite{GL51} of the inverse spectral problem published in 1951 (see also Levitan and Gasymov \cite{LG64}) and alternative variants due to M.~Krein \cite{Kr51}, \cite{Kr53}. For over 45 years, Theorem \ref{t1.1} stood the test of time and resisted any improvements. Finally, in 1998, Simon \cite{Si98} proved the following spectacular result, a local Borg--Marchenko theorem (see part $(i)$ below) and a significant improvement of the original Borg--Marchenko theorem (see part $(ii)$ below):

%%%%%%%%%%%%%%%%%%%%%%%%%%%%%%%%%%%%%
\begin{theorem}\lb{t1.2} ${}$ \\
$(i)$ Let $a>0$, $0<\veps<\pi/2$ and 
suppose that
\begin{equation} \lb{1.4}
|m_1 (z) - m_2(z)| \underset{|z|\to\infty}{=}  O(e^{-2\Im (z^{1/2})a}) 
\end{equation} 
along the ray $\arg(z) = \pi-\veps$. Then
\begin{equation} \lb{1.5}
V_1(x) = V_2 (x) \, \text{ for a.e.\ $x\in [0,a]$.}
\end{equation}
$(ii)$ Let $0<\veps <\pi/2$ and suppose 
that for all $a>0$,
\begin{equation}
|m_1(z) - m_2(z)| \underset{|z|\to\infty}{=} O(e^{-2\Im (z^{1/2})a})   \lb{1.6} 
\end{equation}
along the ray $\arg(z) = \pi -\veps$. Then
\begin{equation}
V_1(x) = V_2(x) \, \text{ for a.e.\  $x\in [0,\infty)$.}    \lb{1.7}
\end{equation}
\end{theorem}
%%%%%%%%%%%%%%%%%%%%%%%%%%%%%%%%%%%%%%
 
The ray $\arg(z) = \pi -\veps$, $0<\veps < \pi/2$ chosen in 
Theorem \ref{t1.2} is of no particular importance. A limit taken along any non-self-intersecting curve $\cC$ going to infinity in the sector 
$\arg(z)\in ((\pi/2)+\veps, \pi -\veps)$ is permissible. For simplicity we only discussed the Dirichlet boundary condition $u(0)=0$ thus far. However, everything extends to the case of general boundary condition $u'(0) + h u(0) = 0$, $h\in\bbR$.  Moreover, the case of a finite interval problem on $[0,b]$, $b\in (0,\infty)$, instead of the half-line $[0,\infty)$ in Theorem \ref{t1.2}\,$(i)$, with $0<a<b$, and a self-adjoint boundary condition at $x=b$ of the type $u'(b) + h_b u(b) = 0$,  $h_b \in \bbR$, can be handled as well. All of this is treated in detail in \cite{GS00a}.

Remarkably enough, the local Borg--Marchenko theorem proven by Simon \cite{Si98} was just a by-product of his new approach to inverse spectral theory for half-line 
Schr\"odinger operators. Actually, Simon's original result in \cite{Si98} was obtained under a bit weaker conditions on $V$; the result as stated in Theorem \ref{t1.2} is taken from \cite{GS00a} (see also \cite{GS00}). While the original proof of the local 
Borg--Marchenko theorem in \cite{Si98} relied on the full power of a new formalism in inverse spectral theory, a short and fairly elementary proof of Theorem \ref{t1.2} was presented in \cite{GS00a}. Without going into further details at this point, we also mention that \cite{GS00a} contains the analog of the local Borg--Marchenko uniqueness result, Theorem \ref{t1.2} for Schr\"odinger operators on the real line. In addition, the case of half-line Jacobi operators and half-line matrix-valued Schr\"odinger operators was dealt with in \cite{GS00a}.

We should also mention some work of Ramm \cite{Ra99}, \cite{Ra00}, who provided a  proof of Theorem~\ref{t1.2}\,$(i)$ under the additional assumption that $V_j$ are 
short-range potentials satisfying $V_j\in L^1([0,\infty); (1+|x|)dx)$, $j=1,2$. 
A very short proof of Theorem \ref{t1.2}, close in spirit to Borg's original paper 
\cite{Bo52}, was subsequently found by Bennewitz \cite{Be01}. Still other proofs were presented in \cite{Ho01} and \cite{Kn01}. Various local and global uniqueness results for matrix-valued Schr\"odinger, Dirac-type, and Jacobi operators were considered in 
\cite{CG02}, \cite{FKRS07}, \cite{GKM02}, \cite{Sa02}, \cite{Sa06}. A local 
Borg--Marchenko theorem for complex-valued potentials has been proved in \cite{BPW02}; the case of semi-infinite Jacobi operators with complex-valued coefficients was studied in  \cite{We04}. This circle of ideas has been reviewed in \cite{Ge07}. 

After this review of Borg--Marchenko-type uniqueness results for Schr\"odinger operators, we now turn to the principal object of our interest in this paper, the so-called CMV operators. CMV operators are a special class of unitary semi-infinite five-diagonal matrices. But for simplicity, we confine ourselves in this introduction to a discussion of CMV operators on $\bbZ$, that is, doubly infinite CMV operators. Denoting by $\D$ the open unit disk in $\bbC$, let $\alpha$ be a sequence of complex numbers in $\D$, $\alpha=\{\al_k\}_{k \in \Z} \subset \D$. The unitary CMV operator $U$ on $\ell^2(\bbZ)$ then can be written as a special five-diagonal doubly infinite matrix in the standard basis of $\ell^2(\bbZ)$ according to \cite[Sects.\ 4.5, 10.5]{Si04}) as in \eqref{B.8}. For the corresponding half-lattice CMV operators $U_{+,k_0}$ 
in $\ell^2([k_0,\infty)\cap\bbZ)$ we refer to \eqref{B.13}--\eqref{B.15}.

The actual history of CMV operators is more involved: The corresponding unitary semi-infinite five-diagonal matrices were first introduced in 1991 by 
Bunse--Gerstner and Elsner \cite{BGE91}, and subsequently treated in detail by Watkins \cite{Wa93} in 1993 (cf.\ the recent discussion in Simon \cite{Si06}). They were subsequently rediscovered by Cantero, Moral, and Vel\'azquez (CMV) in \cite{CMV03}. In \cite[Sects.\ 4.5, 10.5]{Si04}, Simon introduced the corresponding notion of unitary doubly infinite five-diagonal matrices and coined the term ``extended'' CMV matrices. For simplicity, we will just speak of CMV operators whether or not they are half-lattice or full-lattice operators. We also note that in a context different from orthogonal polynomials on the unit circle, Bourget, Howland, and Joye \cite{BHJ03} introduced a family of doubly infinite matrices with three sets of parameters which, for special choices of the parameters, reduces to two-sided CMV matrices on $\bbZ$. Moreover, it is possible to connect unitary block Jacobi matrices to the trigonometric moment problem (and hence to CMV matrices) as discussed by Berezansky and Dudkin \cite{BD05}, \cite{BD06}.

The relevance of this unitary operator $U$ on $\ell^2(\bbZ)$, more precisely, the relevance of the corresponding half-lattice CMV operator $U_{+,0}$ in $\ell^2(\bbN_0)$ is derived from its intimate relationship with the trigonometric moment problem
and hence with finite measures on the unit circle $\dD$. (Here $\bbN_0=\bbN\cup\{0\}$.)  
This will be reviewed in some detail in Section \ref{sB} but we also refer to the
monumental two-volume treatise by Simon \cite{Si04} (see also \cite{Si04b} and 
\cite{Si05}) and the exhaustive bibliography therein. For classical results on orthogonal polynomials on the unit circle we refer, for instance, to \cite{Ak65}, 
\cite{Ge46}--\cite{Ge61}, \cite{Kr45}, \cite{Sz20}--\cite{Sz78}, 
\cite{Ve35}, \cite{Ve36}. More recent references relevant to the spectral theoretic content of this paper are \cite{De07}, \cite{GJ96}--\cite{GT94},  
\cite{GZ06}--\cite{GN01}, \cite{Lu04}, \cite{PY04}, and \cite{Si04a}. The full-lattice CMV operator $U$ on $\bbZ$ is closely related to an important, and only recently intensively studied, completely integrable version of the defocusing nonlinear 
Schr\"odinger equation (continuous in time but discrete in space), a special case of the Ablowitz--Ladik system. Relevant references in this context are, for instance, 
\cite{AL75}--\cite{APT04}, \cite{GGH05}, \cite{GH05}--\cite{GHMT07a}, \cite{Li05}, 
\cite{MEKL95}--\cite{Ne06}, \cite{Sc89}, \cite{Ve99}, and the literature cited therein. 

Next, we briefly summarize some of the principal results proven in this paper. For brevity we just focus on CMV operators on $\bbZ$. We use the following notation for the diagonal and for the neighboring off-diagonal entries of the Green's function of 
$U$ (i.e., the discrete integral kernel of $(U-zI)^{-1}$), 
\begin{align}
g(z,k) &= (U-Iz)^{-1}(k,k),    \quad h(z,k) =
\begin{cases}
(U-Iz)^{-1}(k-1,k), & k \text{ odd}, \\
(U-Iz)^{-1}(k,k-1), & k \text{ even},
\end{cases}\quad k\in\Z,\; z\in\D.     \lb{1.9}
\end{align}

Then the following uniqueness results for CMV operators $U$ on $\bbZ$ will be proven in Section \ref{s2}:

%%%%%%%%%%%%%%%%%%%%%%%%%%%%%%%%%%%%%
\begin{theorem}  \lb{t1.3} 
Assume $\al=\{\al_k\}_{k \in \Z} \subset \D$ and let $k_0\in\Z$. Then any of the
following two sets of data
\begin{enumerate}[$(i)$]
\item $g(z,k_0)$ and $h(z,k_0)$ for all $z$ in a sufficiently small 
neighborhood of the origin under the
assumption $h(0,k_0)\neq0$;
\item $g(z,k_0-1)$ and $g(z,k_0)$ for all $z$ in a sufficiently small 
neighborhood of the origin and $\al_{k_0}$
under the assumption $\al_{k_0}\neq0$;
\end{enumerate}
uniquely determine the Verblunsky coefficients $\{\al_k\}_{k\in\Z}$,
and hence the full-lattice CMV operator $U$.
\end{theorem}
%%%%%%%%%%%%%%%%%%%%%%%%%%%%%%%%%%%%%

In the following local uniqueness result, $g^{(j)}$ and $h^{(j)}$ denote the 
corresponding quantities \eqref{1.9} and \eqref{1.10} associated with the 
Verblunsky coefficients $\alpha^{(j)}$, $j=1,2$. 

%%%%%%%%%%%%%%%%%%%%%%%%%%%%%%%%%%%%%
\begin{theorem}  \lb{t1.4}
Assume $\al^{(\ell)}=\{\al_k^{(\ell)}\}_{k \in \Z} \subset \D$, $\ell=1,2$, 
and let $k_0\in\Z$, $N\in\N$. Then for the full-lattice
problems associated with $\al^{(1)}$ and $\al^{(2)}$, the following
local uniqueness results hold: 
\begin{enumerate}[$(i)$]
\item  If either $h^{(1)}(0,k_0)$ or $h^{(2)}(0,k_0)$ is nonzero and
\begin{align}
\begin{split}
& \abs{g^{(1)}(z,k_0)-g^{(2)}(z,k_0)} +
\abs{h^{(1)}(z,k_0)-h^{(2)}(z,k_0)} \underset{z\to 0}{=} o(z^N),  \\
& \, \text{then } \, \al^{(1)}_k = \al^{(2)}_k \, \text{ for }\, k_0-N \leq k\leq
k_0+N+1.   \lb{1.10}
\end{split}
\end{align}
\item If $\al^{(1)}_{k_0}=\al^{(2)}_{k_0} \neq 0$ and
\begin{align}
\begin{split}
& \abs{g^{(1)}(z,k_0-1)-g^{(2)}(z,k_0-1)} +
\abs{g^{(1)}(z,k_0)-g^{(2)}(z,k_0)} \underset{z\to 0}{=} o(z^N),   \\ 
& \, \text{then } \, \al^{(1)}_k = \al^{(2)}_k \,\text{ for }\, k_0-N-1 \leq k\leq
k_0+N+1.   \lb{1.11}
\end{split}
\end{align}
\end{enumerate}
\end{theorem}
%%%%%%%%%%%%%%%%%%%%%%%%%%%%%%%%%%%%% 

Finally, a brief description of the content of each section in this paper:
In Section \ref{sB} we review the basic Weyl--Titchmarsh theory for CMV operators, discussed in great detail in \cite{GZ06}, as this plays a fundamental role in our principal Section \ref{s2}. In Section \ref{s2} we first provide an alternative proof of the known 
Borg--Marchenko-type uniqueness results for half-lattice CMV operators recorded, for instance, in Simon in \cite[Thm.\ 1.5.5]{Si04} (cf.\ our Theorems \ref{t2.2} and \ref{t2.3}).   
Then we turn to the case of full-lattice CMV operators and prove our principal new results in Theorems \ref{t3.3} and \ref{t3.4} (summarized as 
Theorems \ref{t1.3} and \ref{t1.4} above). In particular, we note that our discussion of CMV operators on the full-lattice will be undertaken in the spirit of \cite{GKM02}, where (local and global) uniqueness theorems for full-line (resp., full-lattice) problems are  formulated in terms of diagonal Green's functions $g(z,x_0)$ and their $x$-derivatives $g^\prime (z,x_0)$ at some fixed $x_0\in\bbR$, for Schr\"odinger and Dirac-type operators on $\bbR$ and similarly for Jacobi operators on $\bbZ$.

An extension of the results of this paper to matrix-valued Verblunsky coefficients appeared in \cite{CGZ07}.

%%%%%%%%%%%%%%%%%%%%%%%%%%%%%%%%%%%%%
%%%%%%%%%%%%%%%%%%%%%%%%%%%%%%%%%%%%%
\section{A Summary of Weyl--Titchmarsh Theory for CMV Operators \\ on
Half-Lattices and on $\bbZ$} \lb{sB}
%%%%%%%%%%%%%%%%%%%%%%%%%%%%%%%%%%%%%%
%%%%%%%%%%%%%%%%%%%%%%%%%%%%%%%%%%%%%%

We start by introducing some of the basic notations used
throughout this paper. Detailed proofs of all facts in this
preparatory section (and a lot of additional material) can be found in
\cite{GZ06}.

In the following, let $\ell^2(\Z)$ be the usual Hilbert space of
all square summable complex-valued sequences with scalar product
$(\cdot,\cdot)_{\ell^2(\bbZ)}$ linear in the second argument.
The {\it standard basis} in $\ell^2(\bbZ)$ is denoted by
\begin{equation}
\{\delta_k\}_{k\in\bbZ}, \quad
\delta_k=(\dots,0,\dots,0,\underbrace{1}_{k},0,\dots,0,\dots)^\top,
\; k\in\bbZ.
\end{equation}
 
For $m\in\N$ and $J\subseteq\bbR$ an interval, we will identify
$\oplus_{j=1}^m\ell^2(J\cap\bbZ)$ and
$\ell^2(J\cap\bbZ)\otimes\bbC^m$ and then use the simplified notation
$\ell^2(J\cap\bbZ)^m$. For simplicity, the identity operators on
$\ell^2(J\cap\bbZ)$ and $\ell^2(J\cap\bbZ)^m$ are abbreviated by $I$
and $I_m$, respectively, without separately indicating its dependence
on $J$. 

By a {\it Laurent polynomial} we denote a finite linear combination
of terms $z^k$, $k\in\bbZ$, with complex-valued coefficients.

Throughout this section we make the following basic assumption:

%%%%%%%%%%%%%%%%%%%%%%%%%%%%%%%%%%%%%%%
\begin{hypothesis} \lb{hB.1}
Let $\alpha$ be a sequence of complex numbers such that
\begin{equation} \lb{B.2}
\alpha=\{\al_k\}_{k \in \Z} \subset \D.
\end{equation}
\end{hypothesis}
%%%%%%%%%%%%%%%%%%%%%%%%%%%%%%%%%%%%%%%

Given a sequence $\alpha$ satisfying \eqref{B.2}, we define the
sequence of positive real numbers $\{\rho_k\}_{k\in\bbZ}$ and
two sequences of complex numbers with positive real parts
$\{a_k\}_{k\in\bbZ}$ and $\{b_k\}_{k\in\bbZ}$ by
\begin{equation}
\rho_k = \sqrt{1-\abs{\al_k}^2}, \quad a_k = 1+\al_k, \quad b_k =
1-\al_k, \quad k \in \Z. \lb{B.3}
\end{equation}

Following Simon \cite{Si04}, we call $\alpha_k$ the Verblunsky
coefficients in honor of Verblunsky's pioneering work in the
theory of orthogonal polynomials on the unit circle \cite{Ve35},
\cite{Ve36}.

Next, we also introduce a sequence of $2\times 2$ unitary
matrices $\Te_k$ by
\begin{equation} \lb{B.4}
\Te_k = \begin{pmatrix} -\al_k & \rho_k \\ \rho_k & \ol{\al_k}
\end{pmatrix},
\quad k \in \Z,
\end{equation}
and two unitary operators $V$ and $W$ on $\ell^2(\Z)$ by their
matrix representations in the standard basis of $\ell^2(\bbZ)$
as follows,
\begin{align} \lb{B.5}
V &= \begin{pmatrix} \ddots & & &
\raisebox{-3mm}[0mm][0mm]{\hspace*{-5mm}\Huge $0$}  \\ &
\Te_{2k-2} & & \\ & & \Te_{2k} & & \\ &
\raisebox{0mm}[0mm][0mm]{\hspace*{-10mm}\Huge $0$} & & \ddots
\end{pmatrix}, \quad
W = \begin{pmatrix} \ddots & & &
\raisebox{-3mm}[0mm][0mm]{\hspace*{-5mm}\Huge $0$}
\\ & \Te_{2k-1} &  &  \\ &  & \Te_{2k+1} &  & \\ &
\raisebox{0mm}[0mm][0mm]{\hspace*{-10mm}\Huge $0$} & & \ddots
\end{pmatrix},
\end{align}
where
\begin{align}
\begin{pmatrix}
V_{2k-1,2k-1} & V_{2k-1,2k} \\
V_{2k,2k-1}   & V_{2k,2k}
\end{pmatrix} =  \Te_{2k},
\quad
\begin{pmatrix}
W_{2k,2k} & W_{2k,2k+1} \\ W_{2k+1,2k}  & W_{2k+1,2k+1}
\end{pmatrix} =  \Te_{2k+1},
\quad k\in\Z.
\end{align}
Moreover, we introduce the unitary operator $U$ on $\ell^2(\Z)$
by
\begin{equation} \lb{B.7}
U = VW,
\end{equation}
or in matrix form, in the standard basis of $\ell^2(\bbZ)$, by
\begin{align}
U &= \begin{pmatrix} \ddots &&\hspace*{-8mm}\ddots
&\hspace*{-10mm}\ddots &\hspace*{-12mm}\ddots
&\hspace*{-14mm}\ddots &&&
\raisebox{-3mm}[0mm][0mm]{\hspace*{-6mm}{\Huge $0$}}
\\
&0& -\al_{0}\rho_{-1} & -\ol{\al_{-1}}\al_{0} & -\al_{1}\rho_{0}
& \rho_{0}\rho_{1}
\\
&& \rho_{-1}\rho_{0} &\ol{\al_{-1}}\rho_{0} &
-\ol{\al_{0}}\al_{1} & \ol{\al_{0}}\rho_{1} & 0
\\
&&&0& -\al_{2}\rho_{1} & -\ol{\al_{1}}\al_{2} & -\al_{3}\rho_{2}
& \rho_{2}\rho_{3}
\\
&&\raisebox{-4mm}[0mm][0mm]{\hspace*{-6mm}{\Huge $0$}} &&
\rho_{1}\rho_{2} & \ol{\al_{1}}\rho_{2} & -\ol{\al_{2}}\al_{3} &
\ol{\al_{2}}\rho_{3}&0
\\
&&&&&\hspace*{-14mm}\ddots &\hspace*{-14mm}\ddots
&\hspace*{-14mm}\ddots &\hspace*{-8mm}\ddots &\ddots
\end{pmatrix}     \lb{B.8}   \\
&= \rho^- \rho \, \deven \, S^{--} + ({\ol \alpha}^-\rho \, \deven - \alpha^+\rho \, \dodd) S^-  - {\ol \alpha}\alpha^+ + ({\ol \alpha} \rho^+ \, \deven - \alpha^{++} \rho^+ \, \dodd) S^+   \no \\
& \quad  + \rho^+ \rho^{++} \, \dodd \, S^{++}, 
\end{align} 
where $\deven$ and $\dodd$ denote the characteristic functions of the even and odd integers,
\begin{equation}
\deven = \chi_{_{2\bbZ}}, \quad \dodd = 1 - \deven = \chi_{_{2\bbZ +1}}.
\end{equation}
Here the diagonal entries in the infinite matrix \eqref{B.8} are given by 
$U_{k,k}=-\ol{\alpha_k}\alpha_{k+1}$, $k\in\Z$.

As explained in the introduction, in the recent literature on orthogonal polynomials on the unit circle, such operators $U$ are frequently called CMV operators.

Next we recall some of the principal results of \cite{GZ06} needed in this paper.

%%%%%%%%%%%%%%%%%%%%%%%%%%%%%%%%%%%%
\begin{lemma} \lb{lB.1}
Let $z\in\bbC\backslash\{0\}$ and $\{u(z,k)\}_{k\in\bbZ},
\{v(z,k)\}_{k\in\bbZ}$ be sequences of complex functions. Then
the following items $(i)$--$(iii)$ are equivalent:
\begin{align}
& (i) \quad (U u(z,\cdot))(k) = z u(z,k), \quad (W u(z,\cdot))=z
v(z,k), \quad k\in\Z.
\\
& (ii) \quad (W u(z,\cdot))(k)=z v(z,k), \quad (V v(z,\cdot))(k)
= u(z,k), \quad k\in\Z.
\\
& (iii) \quad \binom{u(z,k)}{v(z,k)} = T(z,k)
\binom{u(z,k-1)}{v(z,k-1)}, \quad k\in\Z,  \lb{B.11}
\end{align}
where the transfer matrices $T(z,k)$, $z\in\Cz$, $k\in\Z$, are
given by
\begin{equation}
T(z,k) = \begin{cases} \frac{1}{\rho_{k}} \begin{pmatrix}
\al_{k} & z \\ 1/z & \ol{\al_{k}} \end{pmatrix},  & \text{$k$
odd,}  \\ \frac{1}{\rho_{k}} \begin{pmatrix} \ol{\al_{k}} & 1 \\
1 & \al_{k} \end{pmatrix}, & \text{$k$ even.} \end{cases}
\lb{B.12}
\end{equation}
\end{lemma}
%%%%%%%%%%%%%%%%%%%%%%%%%%%%%%%%%%%%

If one sets $\al_{k_0} = e^{is}$, $s\in [0,2\pi)$, for some
reference point $k_0\in\Z$, then the operator $U$ splits into a
direct sum of two half-lattice operators $U_{-,k_0-1}^{(s)}$ and
$U_{+,k_0}^{(s)}$ acting on $\ell^2((-\infty,k_0-1]\cap\Z)$ and
on $\ell^2([k_0,\infty)\cap\Z)$, respectively. Explicitly, one
obtains
\begin{align}
\begin{split}
& U=U_{-,k_0-1}^{(s)} \oplus U_{+,k_0}^{(s)} \, \text{ in } \,
\ell^2((-\infty,k_0-1]\cap\Z) \oplus \ell^2([k_0,\infty)\cap\Z)
\\ & \text{if } \, \al_{k_0} = e^{is}, \; s\in [0,2\pi).
\lb{B.13}
\end{split}
\end{align}
(Strictly, speaking, setting $\al_{k_0} = e^{is}$, $s\in
[0,2\pi)$, for some reference point $k_0\in\Z$ contradicts our
basic Hypothesis \ref{hB.1}. However, as long as the exception
to Hypothesis \ref{hB.1} refers to only one site, we will safely
ignore this inconsistency in favor of the notational simplicity
it provides by avoiding the introduction of a properly modified
hypothesis on $\{\alpha_k\}_{k\in\bbZ}$.) Similarly, one obtains
$W_{-,k_0-1}^{(s)}$, $V_{-,k_0-1}^{(s)}$ and $W_{+,k_0}^{(s)}$,
$V_{+,k_0}^{(s)}$ such that
\begin{equation}
U_{\pm,k_0}^{(s)} = V_{\pm,k_0}^{(s)} W_{\pm,k_0}^{(s)}.
\end{equation}
For simplicity we will abbreviate
\begin{equation}
U_{\pm,k_0} =
U_{\pm,k_0}^{(s=0)}=V_{\pm,k_0}^{(s=0)}W_{\pm,k_0}^{(s=0)}
=V_{\pm,k_0} W_{\pm,k_0}.  \lb{B.15}
\end{equation}

%%%%%%%%%%%%%%%%%%%%%%%%%%%%%%%%%%%%
\begin{lemma} \lb{lB.2}
Let $k_0\in\bbZ$, $z\in\bbC\backslash\{0\}$, and
$\{\hat p_+(z,k,k_0)\}_{k\geq k_0}$, $\{\hat r_+(z,k,k_0)\}_{k\geq k_0}$
be two sequences of complex functions. Then the following items
$(i)$--$(iii)$ are equivalent:
\begin{align}
& (i) \quad (U_{+,k_0} \hat p_+(z,\cdot,k_0))(k) = z \hat p_+(z,k,k_0),
\quad (W_{+,k_0} \hat p_+(z,\cdot,k_0))(k) = z \hat r_+(z,k,k_0), \quad
k\geq k_0.
\\
& (ii) \quad (W_{+,k_0} \hat p_+(z,\cdot,k_0))(k) = z \hat r_+(z,k,k_0),
\quad (V_{+,k_0} \hat r_+(z,\cdot,k_0))(k) = \hat p_+(z,k,k_0), \quad
k\geq k_0.
\\ 
& (iii) \quad \binom{\hat p_+(z,k,k_0)}{\hat r_+(z,k,k_0)} = T(z,k)
\binom{\hat p_+(z,k-1,k_0)}{\hat r_+(z,k-1,k_0)}, \quad k > k_0, \quad
\text{with initial condition}   \no \\
& \hspace{11mm} \hat p_+(z,k_0,k_0) =
\begin{cases}
z \hat r_+(z,k_0,k_0), & \text{$k_0$ odd}, \\
\hat r_+(z,k_0,k_0), & \text{$k_0$ even}.
\end{cases}
\lb{B.18}
\end{align}
Next, consider sequences $\{\hat p_-(z,k,k_0)\}_{k\leq k_0}$,
$\{\hat r_-(z,k,k_0)\}_{k\leq k_0}$. Then the following items
$(iv)$--$(vi)$are equivalent:
\begin{align}
& (iv) \quad (U_{-,k_0} \hat p_-(z,\cdot,k_0))(k) = z \hat p_-(z,k,k_0),
\quad (W_{-,k_0} \hat p_-(z,\cdot,k_0))(k) = z \hat r_-(z,k,k_0), \quad
k\leq k_0.
\\
& (v) \quad  (W_{-,k_0} \hat p_-(z,\cdot,k_0))(k) = z \hat r_-(z,k,k_0),
\quad (V_{-,k_0} \hat r_-(z,\cdot,k_0))(k) = \hat p_-(z,k,k_0), \quad
k\leq k_0.
\\
& (vi) \quad \binom{\hat p_-(z,k-1),k_0}{\hat r_-(z,k-1,k_0)}=T(z,k)^{-1}
\binom{\hat p_-(z,k,k_0)}{\hat r_-(z,k,k_0)}, \quad k \leq k_0, \quad
\text{with initial condition}  \no  
\\
& \hspace*{11mm} \hat p_-(z,k_0,k_0) =
\begin{cases}
\hat r_-(z,k_0,k_0), & \text{$k_0$ odd,}
\\
-z \hat r_-(z,k_0,k_0), & \text{$k_0$ even.}
\end{cases}
 \lb{B.21}
\end{align}
\end{lemma}
%%%%%%%%%%%%%%%%%%%%%%%%%%%%%%%%%%%%

In the following, we denote by
$\Big(\begin{smallmatrix}p_\pm(z,k,k_0)\\
r_\pm(z,k,k_0)\end{smallmatrix}\Big)_{k\in\Z}$ and
$\Big(\begin{smallmatrix}q_\pm(z,k,k_0)\\
s_\pm(z,k,k_0)\end{smallmatrix}\Big)_{k\in\Z}$,
$z\in\bbC\backslash\{0\}$, four linearly independent solutions
of \eqref{B.11} with the following initial conditions:
\begin{align}
\binom{p_+(z,k_0,k_0)}{r_+(z,k_0,k_0)} = \begin{cases}
\binom{z}{1}, &
\text{$k_0$ odd,} \\[1mm]
\binom{1}{1}, & \text{$k_0$ even,} \end{cases} \quad
\binom{q_+(z,k_0,k_0)}{s_+(z,k_0,k_0)} = \begin{cases}
\binom{z}{-1}, &
\text{$k_0$ odd,} \\[1mm]
\binom{-1}{1}, & \text{$k_0$ even.} \end{cases} \lb{B.22}
\\
\binom{p_-(z,k_0,k_0)}{r_-(z,k_0,k_0)} = \begin{cases}
\binom{1}{-1}, &
\text{$k_0$ odd,} \\[1mm]
\binom{-z}{1}, & \text{$k_0$ even,} \end{cases} \quad
\binom{q_-(z,k_0,k_0)}{s_-(z,k_0,k_0)} = \begin{cases}
\binom{1}{1}, &
\text{$k_0$ odd,} \\[1mm]
\binom{z}{1}, & \text{$k_0$ even.} \end{cases} \lb{B.23}
\end{align}
Then it follows that $p_\pm(z,k,k_0)$, $q_\pm(z,k,k_0)$,
$r_\pm(z,k,k_0)$, and $s_\pm(z,k,k_0)$, $k,k_0\in\Z$, are
Laurent polynomials in $z$.

%%%%%%%%%%%%%%%%%%%%%%%%%%%%%%%%%%%%
\begin{lemma} \lb{lB.3}
Let $k_0\in\Z$. Then the sets
$\{p_\pm(\cdot,k,k_0)\}_{k\gtreqless k_0}$ and
$\{r_\pm(\cdot,k,k_0)\}_{k\gtreqless k_0}$ form complete
orthonormal systems of Laurent polynomials in $\Lt{\pm}$, where
\begin{equation}
d\mu_\pm(\zeta,k_0) =
d(\de_{k_0},E_{U_{\pm,k_0}}(\zeta)\de_{k_0})_{\ell^2(\Z\cap[k_0,\pm\infty))}, 
\quad \zeta\in\dD,
\lb{B.24}
\end{equation}
and $dE_{U_{\pm,k_0}}(\cdot)$ denote the operator-valued
spectral measures of the operators $U_{\pm,k_0}$,
\begin{equation}
U_{\pm,k_0}=\oint_{\dD} dE_{U_{\pm,k_0}}(\zeta)\,\zeta.
\end{equation}
Moreover, the half-lattice CMV operators
$U_{\pm,k_0}$ are unitarily equivalent to the operators of
multiplication by the function $id$ $($where $id(\ze)=\ze$,
$\ze\in\dD$$)$ on $\Lt{\pm}$. In particular,
\begin{equation}
\si(U_{\pm,k_0}) =\supp \, (d\mu_\pm(\cdot,k_0))
\end{equation}
and the spectrum of $U_{\pm,k_0}$ is simple. 
\end{lemma}
%%%%%%%%%%%%%%%%%%%%%%%%%%%%%%%%%%%

We note that the measures $d\mu_\pm(\cdot,k_0)$, $k_0\in\bbZ$,
are nonnegative and supported on infinite sets.

%%%%%%%%%%%%%%%%%%%%%%%%%%%%%%%%%%%%
\begin{corollary} \lb{cB.5}
Let $k_0\in\bbZ$. \\ The Laurent polynomials
$\{p_+(\cdot,k,k_0)\}_{k\geq k_0}$ can be constructed by
Gram--Schmidt orthogonalizing
\begin{equation}
\begin{cases} \ze,\,1,\,\ze^2,\,\ze^{-1},\,\ze^3,\,\ze^{-2},\dots, &
\text{$k_0$ odd,} \\
1,\,\ze,\,\ze^{-1},\,\ze^2,\,\ze^{-2},\ze^3,\,\dots, & \text{$k_0$
even}
\end{cases}
\end{equation}
in $\Lt{+}$. \\ The Laurent polynomials $\{r_+(\cdot,k,k_0)\}_{k\geq
k_0}$ can be constructed by Gram--Schmidt orthogonalizing
\begin{equation}
\begin{cases} 1,\,\ze,\,\ze^{-1},\,\ze^2,\,\ze^{-2},\ze^3,\,\dots, &
\text{$k_0$ odd,} \\
1,\,\ze^{-1},\,\ze,\,\ze^{-2},\,\ze^2,\ze^{-3},\,\dots, & \text{$k_0$
even}
\end{cases}
\end{equation}
in $\Lt{+}$. \\ The Laurent polynomials $\{p_-(\cdot,k,k_0)\}_{k\leq
k_0}$ can be constructed by Gram--Schmidt orthogonalizing
\begin{equation}
\begin{cases} 1,\,-\ze,\,\ze^{-1},\,-\ze^2,\,\ze^{-2},-\ze^3,\,\dots, &
\text{$k_0$ odd,} \\
-\ze,\,1,\,-\ze^2,\,\ze^{-1},\,-\ze^3,\,\ze^{-2},\dots, & \text{$k_0$
even}
\end{cases}
\end{equation}
in $\Lt{-}$. \\ The Laurent polynomials $\{r_-(\cdot,k,k_0)\}_{k\leq
k_0}$ can be constructed by Gram--Schmidt orthogonalizing
\begin{equation}
\begin{cases} -1,\,\ze^{-1},\,-\ze,\,\ze^{-2},\,-\ze^2,\ze^{-3},\,\dots, &
\text{$k_0$ odd,} \\
1,\,-\ze,\,\ze^{-1},\,-\ze^2,\,\ze^{-2},-\ze^3,\,\dots, & \text{$k_0$
even}
\end{cases}
\end{equation}
in $\Lt{-}$.
\end{corollary}
%%%%%%%%%%%%%%%%%%%%%%%%%%%%%%%%%%%%

%%%%%%%%%%%%%%%%%%%%%%%%%%%%%%%%%%%%
\begin{theorem} \lb{tB.7}
Let $k_0\in\Z$ and
$d\mu_\pm(\cdot,k_0)$ be nonnegative finite measures on $\dD$
which are supported on infinite sets and normalized by
\begin{align}
\oint_{\dD} d\mu_\pm(\ze,k_0) = 1.
\end{align}
Then $d\mu_\pm(\cdot,k_0)$ are necessarily the spectral measures
for some half-lattice CMV operators $U_{\pm,k_0}$ with
coefficients $\{\al_k\}_{k \geq k_0+1}$, respectively
$\{\al_k\}_{k \leq k_0}$, defined as follows,
\begin{equation}
\al_k = -
\begin{cases}
\big(p_+(\cdot,k-1,k_0),M_{\pm,k_0}(id)
r_+(\cdot,k-1,k_0)\big)_{\Lt{+}}, & k \text{ odd,} \\
\big(r_+(\cdot,k-1,k_0),p_+(\cdot,k-1,k_0)\big)_{\Lt{+}}, & k
\text{ even}
\end{cases}   \lb{B.32}
\end{equation}
for all $k \geq k_0+1$ and
\begin{equation}
\al_k = - \begin{cases} \big(p_-(\cdot,k-1,k_0),M_{\pm,k_0}(id)
r_-(\cdot,k-1,k_0)\big)_{\Lt{-}}, & k \text{ odd,} \\
\big(r_-(\cdot,k-1,k_0),p_-(\cdot,k-1,k_0)\big)_{\Lt{-}}, & k
\text{ even} \end{cases} \lb{B.33}
\end{equation}
for all $k \leq k_0$. Here
$\{p_+(\cdot,k,k_0),r_+(\cdot,k,k_0)\}_{k\geq k_0}$ and
$\{p_-(\cdot,k,k_0),r_-(\cdot,k,k_0)\}_{k\leq k_0}$ denote the
Laurent orthonormal polynomials constructed in Corollary
\ref{cB.5}.
\end{theorem}
%%%%%%%%%%%%%%%%%%%%%%%%%%%%%%%%%%%

Next, we introduce the functions $m_\pm(z,k_0)$ by
\begin{align}
\begin{split}
m_\pm(z,k_0) &= \pm
(\delta_{k_0},(U_{\pm,k_0}+zI)(U_{\pm,k_0}-zI)^{-1}
\delta_{k_0})_{\ell^2(\Z\cap[k_0,\pm\infty))} \\
& =\pm \oint_\dD d\mu_{\pm}(\zeta,k_0)\,
\frac{\zeta+z}{\zeta-z}, \quad z\in\bbC\backslash\dD, 
\end{split} \lb{B.35}
\intertext{with} m_\pm(0,k_0)&=\pm\oint_{\dD} d\mu_\pm(\zeta,k_0)=\pm 1. \lb{B.36}
\end{align}

%%%%%%%%%%%%%%%%%%%%%%%%%%%%%%%%%%%
\begin{theorem} \lb{tB.10}
Let $k_0\in\bbZ$. Then there exist unique functions
$M_\pm(\cdot,k_0)$ such that
\begin{align}
&\binom{u_\pm(z,\cdot,k_0)}{v_\pm(z,\cdot,k_0)} =
\binom{q_+(z,\cdot,k_0)}{s_+(z,\cdot,k_0)} + M_\pm(z,k_0)
\binom{p_+(z,\cdot,k_0)}{r_+(z,\cdot,k_0)} \in
\ell^2([k_0,\pm\infty)\cap\Z)^2, \no \\
& \hspace*{8cm} z\in\bbC\backslash(\dD\cup\{0\}). \lb{B.37}
\end{align}
\end{theorem}
%%%%%%%%%%%%%%%%%%%%%%%%%%%%%%%%%%%

We will call $u_\pm(z,\cdot,k_0)$ (resp., $v_\pm(z,\cdot,k_0)$)
{\it Weyl--Titchmarsh solutions} of $U$ (resp., $U^\top$).
Similarly, we will call $m_\pm(z,k_0)$ as well as $M_\pm(z,k_0)$
the {\it half-lattice Weyl--Titchmarsh $m$-functions} associated
with $U_{\pm,k_0}$. (See also \cite{Si04a} for a comparison of
various alternative notions of Weyl--Titchmarsh $m$-functions
for $U_{+,k_0}$.)

One verifies that
\begin{align}
M_+(z,k_0) &= m_+(z,k_0), \quad z\in\bbC\backslash\dD, \lb{B.38}
\\
M_+(0,k_0) &=1, \lb{B.39}
\\
M_-(z,k_0) &= \frac{\Re(a_{k_0}) +
i\Im(b_{k_0})m_-(z,k_0-1)}{i\Im(a_{k_0}) +
\Re(b_{k_0})m_-(z,k_0-1)} =
\f{(1-z)m_-(z,k_0)+(1+z)}{(1+z)m_-(z,k_0)+(1-z)} \no
\\
&=
\f{(m_-(z,k_0)+1)-z(m_-(z,k_0)-1)}{(m_-(z,k_0)+1)+z(m_-(z,k_0)-1)},
\quad z\in\bbC\backslash\dD, \lb{B.40}
\\
M_-(0,k_0) &=\f{\alpha_{k_0}+1}{\alpha_{k_0}-1}, \lb{B.41}
\\
m_-(z,k_0) &=
\f{\Re(a_{k_0+1})-i\Im(a_{k_0+1})M_-(z,k_0+1)}{\Re(b_{k_0+1})M_-(z,k_0)
- i\Im(b_{k_0+1})} =
\f{(1+z)-(1-z)M_-(z,k_0)}{(1+z)M_-(z,k_0)-(1-z)} \no
\\
&=
\f{z(M_-(z,k_0)+1)-(M_-(z,k_0)-1)}{z(M_-(z,k_0)+1)+(M_-(z,k_0)-1)},
\quad z\in\CdD. \lb{B.42}
\end{align}
In particular, one infers that $M_\pm$ are analytic at $z=0$.

Next, we introduce the functions $\Phi_\pm(\cdot,k)$,
$k\in\bbZ$, by
\begin{align}
\Phi_\pm(z,k) = \f{M_\pm(z,k)-1}{M_\pm(z,k)+1}, \quad
z\in\C\backslash\dD. \lb{B.44}
\end{align}
One then verifies,
\begin{align}
M_\pm(z,k) &= \f{1+\Phi_\pm(z,k)}{1-\Phi_\pm(z,k)}, \quad
z\in\C\backslash\dD, \lb{B.45}
\\
m_-(z,k) &= \f{z-\Phi_-(z,k)}{z+\Phi_-(z,k)}, \quad z\in\CdD.
\lb{B.46}
\end{align}

Finally, we turn to the resolvent of $U$:

%%%%%%%%%%%%%%%%%%%%%%%%%%%%%%%%%%%%
\begin{lemma}
Let $z\in\bbC\backslash(\dD\cup\{0\})$ and fix $k_0\in\bbZ$. Then  the 
resolvent $(U-zI)^{-1}$ of the unitary
CMV operator $U$ on $\ell^2(\bbZ)$ is given in terms of its
matrix representation in the standard basis of $\ell^2(\Z)$ by
\begin{align} \nonumber
(U-zI)^{-1}(k,k') &= \frac{-1}{2z[M_+(z,k_0)-M_-(z,k_0)]} \no \\
& \quad \times \begin{cases} \wti u_-(z,k,k_0)v_+(z,k',k_0), & k
< k' \text{ and } k = k'
\text{ odd}, \\
v_-(z,k',k_0) \wti u_+(z,k,k_0), & k' < k \text{ and } k = k'
\text{ even},
\end{cases} \quad k,k' \in\Z, \lb{B.84}
\end{align}
where 
\begin{align}
\wti u_\pm(z,k,k_0) &=
\begin{cases}
u_\pm(z,k,k_0)/z, & \text{$k_0$ odd,}
\\
u_\pm(z,k,k_0), &  \text{$k_0$ even.}
\end{cases}
\end{align}
Moreover, since $0\in\bbC\backslash\sigma(U)$, \eqref{B.84}
analytically extends to $z=0$. In addition, the following formulas
hold,
\begin{align}
(U-zI)^{-1}(k,k) &=
\frac{[1-M_+(z,k)][1+M_-(z,k)]}{2z[M_+(z,k)-M_-(z,k)]}, \lb{B.86}
\\
(U-zI)^{-1}(k-1,k-1) &= \f{[\ol{a_{k}}-\ol{b_{k}}M_+(z,k)][a_{k}
+b_{k}M_-(z,k)]}{2z\rho_{k}^2[M_+(z,k)-M_-(z,k)]}, \lb{B.86a}
\\
(U-zI)^{-1}(k-1,k) &= - \f{
\begin{cases}
{[1-M_+(z,k)][\ol{a_{k}}-\ol{b_{k}}M_-(z,k)]}, & \text{$k$ odd,}
\\
{[1+M_+(z,k)][a_{k}+b_{k}M_-(z,k)]}, & \text{$k$ even,}
\end{cases}}{2z\rho_{k}[M_+(z,k)-M_-(z,k)]} \lb{B.86b}
\\
(U-zI)^{-1}(k,k-1) &= - \f{
\begin{cases}
{[1+M_+(z,k)][a_{k}+b_{k}M_-(z,k)]}, & \text{$k$ odd,}
\\
{[1-M_+(z,k)][\ol{a_{k}}-\ol{b_{k}}M_-(z,k)]}, & \text{$k$ even.}
\end{cases}}{2z\rho_{k}[M_+(z,k)-M_-(z,k)]} \lb{B.86c}
\end{align}
\end{lemma}
%%%%%%%%%%%%%%%%%%%%%%%%%%%%%%%%%%%

%%%%%%%%%%%%%%%%%%%%%%%%%%%%%%%%%%%%%%%%
%%%%%%%%%%%%%%%%%%%%%%%%%%%%%%%%%%%%%%%%
\section{Borg--Marchenko-type Uniqueness Results for CMV Operators} \lb{s2}
%%%%%%%%%%%%%%%%%%%%%%%%%%%%%%%%%%%%%%%%
%%%%%%%%%%%%%%%%%%%%%%%%%%%%%%%%%%%%%%%%

In this section we prove (local) Borg--Marchenko-type uniqueness results for
CMV operators with scalar-valued Verblunsky coefficients on the full lattice 
$\ell^2(\Z)$ and on half-lattices $\ell^2([k_0,\pm\infty)\cap\Z)$. The principal results in the full lattice case, Theorems \ref{t3.3} and \ref{t3.4} are new. We freely use the notation established in Section \ref{sB}.

We start with uniqueness results for CMV operators on half-lattices. While these results are known and have recently been recorded by Simon in \cite[Thm.\ 1.5.5]{Si04}, we present the proofs for the convenience of the reader as the half-lattice results are crucial ingredients for our new full lattice results. 

%%%%%%%%%%%%%%%%%%%%%%%%%%%%%%%%%%%%%%
\begin{theorem} \lb{t2.2}
Assume Hypothesis \ref{hB.1} and let $k_0\in\Z$, $N\in\N$. Then,
for the right half-lattice problem, the following sets of data $(i)$--$(v)$  
are equivalent:
\begin{align}
(i) &\quad \big\{\al_{k_0+k}\big\}_{k=1}^N.
\\
(ii) &\quad \bigg\{\oint_\dD \ze^k\,d\mu_+(\ze,k_0)\bigg\}_{k=1}^N.
\\
(iii) &\quad \big\{m_{+,k}(k_0)\big\}_{k=1}^N, \;\text{ where
$m_{+,k}(k_0)$, $k\geq 0$, are the Taylor coefficients of
$m_+(z,k_0)$} \no\\ &\hspace*{5mm} \text{at $z=0$, that is, }\;
m_+(z,k_0) = \sum\nolimits_{k=0}^\infty m_{+,k}(k_0)z^k, \; z\in\D.
\\
(iv) &\quad \big\{M_{+,k}(k_0)\big\}_{k=1}^N, \;\text{ where
$M_{+,k}(k_0)$, $k\geq 0$, are the Taylor coefficients of
$M_+(z,k_0)$} \no\\ &\hspace*{5mm} \text{at $z=0$, that is, }\;
M_+(z,k_0) = \sum\nolimits_{k=0}^\infty M_{+,k}(k_0)z^k, \; z\in\D.
\\
(v) &\quad \big\{\phi_{+,k}(k_0)\big\}_{k=1}^N, \;\text{ where
$\phi_{+,k}(k_0)$, $k\geq 0$, are the Taylor coefficients of
$\Phi_+(z,k_0)$} \no\\ &\hspace*{5mm} \text{at $z=0$, that is, }\;
\Phi_+(z,k_0) = \sum\nolimits_{k=0}^\infty \phi_{+,k}(k_0)z^k, \;
z\in\D.
\end{align}
Similarly, for the left half-lattice problem, the following sets
of data $(vi)$--$(x)$ are equivalent:
\begin{align}
(vi) &\quad \big\{\al_{k_0-k}\big\}_{k=0}^{N-1}.
\\
(vii) &\quad \bigg\{\oint_\dD \ze^k\,d\mu_-(\ze,k_0)\bigg\}_{k=1}^N. 
\\
(viii) &\quad \big\{m_{-,k}(k_0)\big\}_{k=1}^{N}, \;\text{ where
$m_{-,k}(k_0)$, $k\geq 0$, are the Taylor coefficients of
$m_-(z,k_0)$} \no\\ &\hspace*{5mm} \text{at $z=0$, that is, }\;
m_-(z,k_0) = \sum\nolimits_{k=0}^\infty m_{-,k}(k_0)z^k. 
\\
(ix) &\quad \big\{M_{-,k}(k_0)\big\}_{k=0}^{N-1}, \;\text{ where
$M_{-,k}(k_0)$, $k\geq 0$, are the Taylor coefficients of
$M_-(z,k_0)$} \no\\ &\hspace*{5mm} \text{at $z=0$, that is, }\;
M_-(z,k_0) = \sum\nolimits_{k=0}^\infty M_{-,k}(k_0)z^k. 
\\
(x) &\quad \big\{\varphi_{-,k}(k_0)\big\}_{k=0}^{N-1}, \;\text{ where
$\varphi_{-,k}(k_0)$, $k\geq 0$, are the Taylor coefficients of
$\Phi_-(z,k_0)^{-1}$} \no\\ &\hspace*{5mm} \text{at $z=0$, that is, }\;
\Phi_-(z,k_0)^{-1} = \sum\nolimits_{k=0}^\infty \varphi_{-,k}(k_0) z^k.
\end{align}
\end{theorem}
%%%%%%%%%%%%%%%%%%%%%%%%%%%%%%%%%%%%%%
\begin{proof}
The crucial equivalence of items $(i)$ and $(ii)$ can be found in 
Simon \cite[Thm.\ 1.5.5]{Si04} (where a more general result is proven). For the convenience of the reader we present an alternative proof below. \\ 
$(i)\Rightarrow(ii)$ and $(vi)\Rightarrow(vii)$: First, utilizing 
relations \eqref{B.18} and \eqref{B.21} with the initial conditions
\eqref{B.22} and \eqref{B.23}, one constructs $\{p_\pm(z,k_0\pm
k,k_0)\}_{k=1}^{N}$ and $\{r_\pm(z,k_0\pm k,k_0)\big\}_{k=1}^{N}$.
We note that the polynomials
\begin{align}
&\begin{cases}%
z^{-1} p_+(z,k_0+k,k_0), \; r_-(z,k_0-k,k_0), & k_0 \text{ odd},
\\[1mm]
r_+(z,k_0+k,k_0), \; z^{-1} p_-(z,k_0-k,k_0), & k_0 \text{ even},
\end{cases} \lb{2.15}
\intertext{are linear combinations of }
&\begin{cases}%
1,z^{-1},z,z^{-2},z^2,\dots,z^{(k-1)/2},z^{-(k+1)/2}, & k \text{
odd},
\\[1mm]
1,z^{-1},z,z^{-2},z^2,\dots,z^{-k/2},z^{k/2}, & k \text{ even},
\end{cases} \lb{2.16}
\intertext{and}
&\begin{cases}%
r_+(z,k_0+k,k_0), \; p_-(z,k_0-k,k_0), & k_0 \text{ odd},
\\[1mm]
p_+(z,k_0+k,k_0), \; r_-(z,k_0-k,k_0), & k_0 \text{ even},
\end{cases} \lb{2.17}
\intertext{are linear combinations of }
&\begin{cases}%
1,z,z^{-1},z^2,z^{-2},\dots,z^{-(k-1)/2},z^{(k+1)/2}, & k \text{ odd},
\\[1mm]
1,z,z^{-1},z^2,z^{-2},\dots,z^{k/2},z^{-k/2}, & k \text{ even}.
\end{cases} \lb{2.18}
\end{align}
Moreover, the last elements of the sequences in \eqref{2.16} and
\eqref{2.18} represent the leading-order terms of the polynomials in
\eqref{2.15} and \eqref{2.17}, respectively, and the corresponding
leading-order coefficients are nonzero.

Next, assume $k_0$ and $k$ to be odd. Then utilizing \eqref{2.17}
and \eqref{2.18} one finds constants $c_{\pm,j}$ and $d_{\pm,j}$,
$0\leq j\leq k$, such that
\begin{align}
z^{-(k-1)/2}  &= \sum_{j=0}^{k} c_{+,j}\, r_+(z,k_0+j, k_0), \quad
z^{(k+1)/2} = \sum_{j=0}^{k} d_{+,j}\, r_+(z,k_0+j,k_0),
\\
z^{-(k-1)/2}  &= \sum_{j=0}^{k} c_{-,j}\, p_-(z,k_0-j, k_0), \quad
z^{(k+1)/2} = \sum_{j=0}^{k} d_{-,j}\, p_-(z,k_0-j,k_0),
\end{align}
and, using  Lemma \ref{lB.3}, computes
\begin{align}
\oint_\dD \ze^k d\mu_\pm(\ze,k_0) = \oint_\dD \ol{\ze^{-(k-1)/2}}
\ze^{(k+1)/2} d\mu_\pm(\ze,k_0) = \sum_{j=0}^k
\ol{c_{\pm,j}}\,d_{\pm,j}.
\end{align}
The remaining cases of $k_0$ and $k$ follow similarly.

$(ii)\Rightarrow(i)$ and $(vii)\Rightarrow(vi)$: Since the measures
$d\mu_\pm(\cdot,k_0)$ are real-valued and normalized, one has
\begin{align}
\oint_\dD \ze^{-k} d\mu_\pm(\ze,k_0) = \ol{\oint_\dD \ze^{k}
d\mu_\pm(\ze,k_0)} \quad\text{and}\quad \oint_\dD d\mu_\pm(\ze,k_0) =
1,
\end{align}
that is, the knowledge of positive moments imply the knowledge of
negative ones. Applying Corollary \ref{cB.5} and Theorem
\ref{tB.7} one constructs orthonormal polynomials $\{p_\pm(\ze,k_0\pm
k,k_0)\}_{k=1}^{N}$ and $\{r_\pm(\ze,k_0\pm k,k_0)\big\}_{k=1}^{N}$
and subsequently the Verblunsky coefficients in $(i)$ and $(vi)$ using formulas
\eqref{B.32} and \eqref{B.33}.

$(ii)\Leftrightarrow(iii)$ and $(vii)\Leftrightarrow(viii)$: These equivalences 
follow directly from \eqref{B.35},
\begin{align}
m_\pm(z,k_0) &= \pm \oint_\dD \frac{\zeta+z}{\zeta-z}\,
d\mu_{\pm}(\zeta,k_0) = \pm 1 \pm 2\sum_{k=1}^{\infty} z^k
\ol{\oint_\dD \ze^{k} d\mu_\pm(\ze,k_0)}, \quad z\in\D.
\end{align}

$(iii)\Leftrightarrow(iv)$: This follows from \eqref{B.38}.

$(iv)\Leftrightarrow(v)$: This follows from \eqref{B.39},
\eqref{B.44}, \eqref{B.45}, and the fact that $|\Phi_+(z,k_0)|<1$,
$z\in\D$,
\begin{align}
M_+(z,k_0) &= \f{1+\Phi_+(z,k_0)}{1-\Phi_+(z,k_0)} \underset{z\to 0}{=}
[1+\Phi_+(z,k_0)] \sum_{k=0}^\infty \Phi_+(z,k_0)^k, 
\\
\Phi_+(z,k_0) &= \f{2^{-1} [M_+(z,k_0)-1]}{2^{-1}[M_+(z,k_0)-1]+1} \underset{z\to 0}{=} 
-\sum_{k=1}^\infty 2^{-k} [1-M_+(z,k_0)]^k.
\end{align}

$(ix)\Leftrightarrow(x)$: This follows from \eqref{B.41},
\eqref{B.44}, \eqref{B.45}, and the fact that
$|\Phi_-(z,k_0)^{-1}|<1$, $z\in\D$,
\begin{align}
M_-(z,k_0) &= \f{1/\Phi_-(z,k_0)+1}{1/\Phi_-(z,k_0)-1} \underset{z\to 0}{=} 
- [1/\Phi_-(z,k_0)+1] \sum_{k=0}^\infty \Phi_-(z,k_0)^{-k},
\\
1/\Phi_-(z,k_0) &=
\f{M_-(z,k_0)+1}{M_-(z,k_0)-M_-(0,k_0)+M_-(0,k_0)-1} \no
\\ &=
\f{[M_-(z,k_0)+1][M_-(0,k_0)-1]^{-1}}
{[M_-(z,k_0)-M_-(0,k_0)][M_-(0,k_0)-1]^{-1}+1}
\\ &  \hspace{-1.5mm} 
\underset{z\to 0}{=} 
\f{M_-(z,k_0)+1}{M_-(0,k_0)-1}\sum_{k=0}^\infty
\bigg(\f{M_-(z,k_0)-M_-(0,k_0)}{1-M_-(0,k_0)}\bigg)^k. \no
\end{align}

$(viii)\Leftrightarrow(x)$: This follows from \eqref{B.36},
\eqref{B.46}, and the fact that $|\Phi_-(z,k_0)^{-1}|<1$, $z\in\D$,
\begin{align}
m_-(z,k_0) &= \f{z/\Phi_-(z,k_0)-1}{z/\Phi_-(z,k_0)+1} \underset{z\to 0}{=} 
[z/\Phi_-(z,k_0)-1] \sum_{k=0}^\infty \big(-z/\Phi_-(z,k_0)\big)^k,
\\
z/\Phi_-(z,k_0) &= \f{1+m_-(z,k_0)}{1-m_-(z,k_0)} \underset{z\to 0}{=} 
\f{2^{-1} [1+m_-(z,k_0)]}{1- 2^{-1} [1+m_-(z,k_0)]} = \sum_{k=1}^\infty
\bigg(\f{1+m_-(z,k_0)}{2}\bigg)^k.
\end{align}

\end{proof}
%%%%%%%%%%%%%%%%%%%%%%%%%%%%%%%%%%%%%%%

We restate Theorem \ref{t2.2} as follows:

%%%%%%%%%%%%%%%%%%%%%%%%%%%%%%%%%%%%%%%
\begin{theorem} \lb{t2.3}
Assume Hypothesis \ref{hB.1} for two sequences $\al^{(1)}$,
$\al^{(2)}$ and let $k_0\in\Z$, $N\in\N$. Then for the right
half-lattice problems associated with $\al^{(1)}$ and $\al^{(2)}$ 
the following items $(i)$--$(iv)$ are equivalent:
\begin{align}
(i) &\quad \al_k^{(1)} = \al_k^{(2)}, \quad k_0+1\leq k\leq k_0+N.
\\
(ii) &\quad m_+^{(1)}(z,k_0)-m_+^{(2)}(z,k_0) \underset{z\to 0}{=} \oh(z^N).
\\
(iii) &\quad M_+^{(1)}(z,k_0)-M_+^{(2)}(z,k_0) \underset{z\to 0}{=} \oh(z^N). 
\\
(iv) &\quad \Phi_+^{(1)}(z,k_0)-\Phi_+^{(2)}(z,k_0) \underset{z\to 0}{=} \oh(z^N).
\end{align}
Similarly, for the left half-lattice problems associated with
$\al^{(1)}$ and $\al^{(2)}$, the following items $(v)$--$(viii)$ are equivalent:
\begin{align}
(v) &\quad \al_k^{(1)} = \al_k^{(2)}, \quad k_0-N+1\leq k\leq k_0. 
\\
(vi) &\quad m_-^{(1)}(z,k_0)-m_-^{(2)}(z,k_0) \underset{z\to 0}{=} \oh(z^N). 
\\
(vii) &\quad M_-^{(1)}(z,k_0)-M_-^{(2)}(z,k_0) \underset{z\to 0}{=} \oh(z^{N-1}). 
\\
(viii) &\quad 1/\Phi_-^{(1)}(z,k_0)-1/\Phi_-^{(2)}(z,k_0) \underset{z\to 0}{=} 
\oh(z^{N-1}).
\end{align}
\end{theorem}
%%%%%%%%%%%%%%%%%%%%%%%%%%%%%%%%%%%%%
\begin{proof}
This follows immediately from Theorem \ref{t2.2}.
\end{proof}
%%%%%%%%%%%%%%%%%%%%%%%%%%%%%%%%%%%%%

Finally, we turn to CMV operators on $\bbZ$.

To start, we introduce the following notation for
the diagonal and for the neighboring off-diagonal entries of the
Green's function of $U$ (i.e., the discrete integral kernel of $(U-zI)^{-1}$),
\begin{align}
g(z,k) &= (U-Iz)^{-1}(k,k),     \lb{2.36}
\\
h(z,k) &=
\begin{cases}
(U-Iz)^{-1}(k-1,k), & k \text{ odd}, \\
(U-Iz)^{-1}(k,k-1), & k \text{ even},
\end{cases}\quad k\in\Z,\; z\in\D.     \lb{2.37}
\end{align}
Then the following uniqueness results hold for the full-lattice CMV
operator $U$: 

%%%%%%%%%%%%%%%%%%%%%%%%%%%%%%%%%%%%%
\begin{theorem}  \lb{t3.3}
Assume Hypothesis \ref{hB.1} and let $k_0\in\Z$. Then any of the
following two sets of data
\begin{enumerate}[$(i)$]
\item $g(z,k_0)$ and $h(z,k_0)$ for all $z$ in some open $($nonempty$)$  
neighborhood of the origin under the
assumption $h(0,k_0)\neq0$;
\item $g(z,k_0-1)$ and $g(z,k_0)$ for all $z$ in some open $($nonempty$)$  
neighborhood of the origin and $\al_{k_0}$
under the assumption $\al_{k_0}\neq0$;
\end{enumerate}
uniquely determines the Verblunsky coefficients $\{\al_k\}_{k\in\Z}$,
and hence the full-lattice CMV operator $U$.
\end{theorem}
%%%%%%%%%%%%%%%%%%%%%%%%%%%%%%%%%%%%%
\begin{proof}
{\it Case} (i). First, note that it follows from \eqref{B.8} that
\begin{align}
g(0,k_0) &= (U^{-1})_{k_0,k_0} = (U^*)_{k_0,k_0} = \ol{U_{k_0,k_0}}
= -\al_{k_0}\ol{\al_{k_0+1}}, \lb{2.38}
\\
h(0,k_0) &= \begin{cases} (U^{-1})_{k_0-1,k_0} = \ol{U_{k_0,k_0-1}} =
-\ol{\al_{k_0+1}}\rho_{k_0}, & k_0 \text{ odd},  
\\
(U^{-1})_{k_0,k_0-1} = \ol{U_{k_0-1,k_0}} =
-\ol{\al_{k_0+1}}\rho_{k_0}, & k_0 \text{ even}. \end{cases}  \lb{2.40}
\end{align}
Since by hypothesis $h(0,k_0)\neq 0$, one can solve the above equalities for
$\al_{k_0}$, 
\begin{align}
\f{g(0,k_0)}{h(0,k_0)} = \al_{k_0}/\rho_{k_0}, \quad \abs{\al_{k_0}}^2 =
\f{\abs{g(0,k_0)}^2}{\abs{g(0,k_0)}^2+\abs{h(0,k_0)}^2}, \quad
\rho_{k_0} =
\f{\abs{h(0,k_0)}}{\sqrt{\abs{g(0,k_0)}^2+\abs{h(0,k_0)}^2}},
\end{align}
and hence,
\begin{align}
\al_{k_0} = \f{g(0,k_0)\abs{h(0,k_0)}}{h(0,k_0)
\sqrt{\abs{g(0,k_0)}^2+\abs{h(0,k_0)}^2}}.
\end{align}
Recalling \eqref{B.3}, one has $a_{k_0}=1+\al_{k_0}$ and
$b_{k_0}=1-\al_{k_0}$.

Next, utilizing \eqref{B.86}, \eqref{B.86b}, and \eqref{B.86c}, one
computes,
\begin{align}
\f{g(z,k_0)}{h(z,k_0)} =
\f{\rho_{k_0}[1+M_-(z,k_0)]}{\ol{b_{k_0}}M_-(z,k_0)-\ol{a_{k_0}}},
\quad z\in\D.
\end{align}
Solving for $M_-(z,k_0)$, one then obtains
\begin{align}
M_-(z,k_0) = \f{2g(z,k_0)}{\ol{b_{k_0}} g(z,k_0) - \rho_{k_0}
h(z,k_0)}-1, \quad z\in\D. \lb{2.44}
\end{align}
Here, the denominator may have only a discrete set of zeros 
(i.e., without accumulation points in $\D$), 
corresponding to removable singularities of the fraction. Otherwise, 
the denominator is identically zero
in $\D$ since the functions $g$ and $h$ are analytic in $\D$. This in turn 
implies that the numerator is identically zero in $\D$, and hence,
$h$ is identically zero in $\D$, contradicting our assumption
$h(0,k_0)\neq 0$.

Next, having obtained $M_-(z,k_0)$, one solves
\begin{align}
h(z,k_0) = - \f{[1-M_+(z,k_0)][\ol{a_{k}}-\ol{b_{k}}M_-(z,k_0)]}
{2z\rho_{k}[M_+(z,k_0)-M_-(z,k_0)]}, \quad z\in\D
\end{align}
for $M_+(z,k_0)$ and obtains,
\begin{align}
M_+(z,k_0) &= \f{[\ol{a_{k_0}}-\ol{b_{k_0}}M_-(z,k_0)] -
2z\rho_{k_0} h(z,k_0) M_-(z,k_0)}
{[\ol{a_{k_0}}-\ol{b_{k_0}}M_-(z,k_0)] - 2z\rho_{k_0} h(z,k_0)} \no
\\
&= \f{2(1+zg(z,k_0))}{1+z [\ol{b_{k_0}}g(z,k_0) - \rho_{k_0}
h(z,k_0)]} -1, \quad z\in\D. \lb{2.46}
\end{align}
Here, the denominator may have only a discrete set of zeros,  
corresponding to removable singularities of the fraction \eqref{2.46}. 
Otherwise, one concludes again that the denominator is identically zero
in $\D$, contradicting the fact that $g$ and $h$ are analytic in $\D$.

Finally, Theorem \ref{t2.2} (parts $(i)$, $(iv)$ and $(vi)$, $(ix)$)
implies that $M_\pm(z,k_0)$, $z\in\D$, uniquely determine the 
Verblunsky coefficients $\{\al_k\}_{k\in\Z}$.

{\it Case} (ii). First, using \eqref{B.44} and
\begin{align}
1+zg(z,k_0) =
\frac{[1+M_+(z,k_0)][1-M_-(z,k_0)]}{2[M_+(z,k_0)-M_-(z,k_0)]},
\end{align}
which follows from \eqref{B.86}, one rewrites \eqref{B.86} and
\eqref{B.86a} as
\begin{align}
zg(z,k_0) &= [\Phi_+(z,k_0)/\Phi_-(z,k_0)] [1+zg(z,k_0)]. \lb{2.49}  \\  
z\rho_{k_0}^2g(z,k_0-1) &=
\frac{\big[[1-M_+(z,k_0)]+\ol{\al_{k_0}}[1+M_+(z,k_0)]\big]
\big[[1+M_-(z,k_0)]+\al_{k_0}[1-M_-(z,k_0)]\big]}
{2[M_+(z,k_0)-M_-(z,k_0)]} \no
\\
&= [\ol{\al_{k_0}}-\Phi_+(z,k_0)][\al_{k_0}-1/\Phi_-(z,k_0)]
[1+zg(z,k_0)], \lb{2.48}
\end{align}

Next, introducing the functions $A(z,k_0)$ and $B(z,k_0)$ by
\begin{align}
A(z,k_0) &= 1+zg(z,k_0), \quad B(z,k_0) =
z\rho_{k_0}^2g(z,k_0-1)-zg(z,k_0)-\abs{\al_{k_0}}^2A(z,k_0),
\end{align}
one obtains from \eqref{2.48} and \eqref{2.49} that $\Phi_+(z,k_0)$
satisfies the following quadratic equation
\begin{align}
\al_{k_0}A(z,k_0)\Phi_+(z,k_0)^2 + B(z,k_0)\Phi_+(z,k_0) +
\ol{\al_{k_0}}zg(z,k_0) = 0. \lb{2.51}
\end{align}
In addition, it follows from \eqref{B.44} that $\Phi_+(0,k_0)=0$
since by \eqref{B.39} $M_+(0,k_0)=1$. Hence, $\Phi_+(z,k_0)$ can be
uniquely determined for sufficiently small $\abs{z}$ from
\eqref{2.51},
\begin{align}
\Phi_+(z,k_0) = \frac{- B(z,k_0) -
\sqrt{B(z,k_0)^2 - 4 \abs{\al_{k_0}}^2A(z,k_0)zg(z,k_0)}}
{2\al_{k_0}A(z,k_0)}. \lb{2.52}
\end{align}
Utilizing \eqref{2.48}, one also finds $1/\Phi_-(z,k_0)$ for
sufficiently small $\abs{z}$,
\begin{align}
1/\Phi_-(z,k_0) = \al_{k_0} - \frac{z\rho_{k_0}^2g(z,k_0-1)}
{[1+zg(z,k_0)][\ol{\al_{k_0}}-\Phi_+(z,k_0)]}. \lb{2.53}
\end{align}

Finally, Theorem \ref{t2.2} (parts $(i)$, $(v)$ and $(vi)$, $(x)$) implies again that for 
$\abs{z}$ sufficiently small, $\Phi_\pm(z,k_0)^{\pm1}$ uniquely determine the Verblunsky coefficients $\{\al_k\}_{k\in\Z}$.
\end{proof}
%%%%%%%%%%%%%%%%%%%%%%%%%%%%%%%%%%%%%

In the following result, $g^{(j)}$ and $h^{(j)}$ denote the corresponding quantities 
\eqref{2.36} and \eqref{2.37} associated with the Verblunsky coefficients 
$\alpha^{(j)}$, $j=1,2$. 

%%%%%%%%%%%%%%%%%%%%%%%%%%%%%%%%%%%%%
\begin{theorem}  \lb{t3.4}
Assume Hypothesis \ref{hB.1} for two sequences $\al^{(1)}$,
$\al^{(2)}$ and let $k_0\in\Z$, $N\in\N$. Then for the full-lattice
problems associated with $\al^{(1)}$ and $\al^{(2)}$, the following
local uniqueness results hold: 
\begin{enumerate}[$(i)$]
\item  If either $h^{(1)}(0,k_0)$ or $h^{(2)}(0,k_0)$ is nonzero and
\begin{align}
\begin{split}
& \abs{g^{(1)}(z,k_0)-g^{(2)}(z,k_0)} +
\abs{h^{(1)}(z,k_0)-h^{(2)}(z,k_0)} \underset{z\to 0}{=} o(z^N), \lb{2.54} \\
& \, \text{then } \, \al^{(1)}_k = \al^{(2)}_k \, \text{ for }\, k_0-N \leq k\leq
k_0+N+1.
\end{split}
\end{align}
\item If $\al^{(1)}_{k_0}=\al^{(2)}_{k_0} \neq 0$ and
\begin{align}
\begin{split}
& \abs{g^{(1)}(z,k_0-1)-g^{(2)}(z,k_0-1)} +
\abs{g^{(1)}(z,k_0)-g^{(2)}(z,k_0)} \underset{z\to 0}{=} o(z^N), \lb{2.56}  \\
& \, \text{then } \, \al^{(1)}_k = \al^{(2)}_k \quad\text{for}\quad k_0-N-1 \leq k\leq
k_0+N+1.
\end{split}
\end{align}
\end{enumerate}
\end{theorem}
%%%%%%%%%%%%%%%%%%%%%%%%%%%%%%%%%%%%% 
\begin{proof}
{\it Case} $(i)$. The result follows from Theorem \ref{t2.3} (parts
$(i)$, $(iii)$ and $(v)$, $(vii)$) upon verifying that \eqref{2.44},
\eqref{2.46}, and \eqref{2.54} imply
\begin{align}
M_+^{(1)}(z,k_0)-M_+^{(2)}(z,k_0) &\underset{z\to 0}{=} o(z^{N+1}),
\\
M_-^{(1)}(z,k_0)-M_-^{(2)}(z,k_0) &\underset{z\to 0}{=} o(z^N).
\end{align}
The latter asymptotic behavior follows from the fact that the denominator in
\eqref{2.54} is non-zero as $z\to0$. Indeed, using
\eqref{2.38} and \eqref{2.40} one computes
\begin{align}
\ol{b_{k_0}}g(0,k_0) - \rho_{k_0}h(0,k_0) =
(\ol{\al_{k_0}}-1)\al_{k_0}\ol{\al_{k_0+1}} +
\rho_{k_0}^2\ol{\al_{k_0+1}} = (\al_{k_0}-1)h(0,k_0)/\rho_{k_0}
\neq0.
\end{align}

{\it Case} $(ii)$. The result follows from Theorem \ref{t2.3} (parts
$(i)$, $(iv)$ and $(v)$, $(viii)$) upon verifying that \eqref{2.52},
\eqref{2.53}, and \eqref{2.56} imply
\begin{align}
\abs{\Phi_+^{(1)}(z,k_0)-\Phi_+^{(2)}(z,k_0)} +
\abs{1/\Phi_-^{(1)}(z,k_0)-1/\Phi_-^{(2)}(z,k_0)} \underset{z\to 0}{=} o(z^{N+1}).
\end{align}
\end{proof}
%%%%%%%%%%%%%%%%%%%%%%%%%%%%%%%%%%%%%%

%%%%%%%%%%%%%%%%%%%%%%%%%%%%%%%%%%%%
%%%%%%%%%%%%%%%%%%%%%%%%%%%%%%%%%%%%

\end{document}